\documentclass{article}
\usepackage{amssymb,amsmath,amsthm}
\usepackage{bm}
\usepackage{graphicx}
\font\cyr=wncyr8

\def\fsf{\{4,6\,|\,4\}}
\def\sff{\{6,4\,|\,4\}}

\def\bR{{\mathbb R}}

\def\bT{{\mathbb T}}
\def\bZ{{\mathbb Z}}
\def\cA{{\cal A}}
\def\cC{{\cal C}}
\def\cD{{\cal D}}
\def\zone{\cD}
\def\cE{{\cal E}}
\def\cM{M}

\def\cP{{\cal P}}

\def\cT{{\cal T}}

\def\mf{\theta}

\def\Tt{\bT^3}

\def\RPt{\bR\hbox{\rm P}^2}

\def\rmU{\uppercase\expandafter{\romannumeral1}}
\def\rmD{\uppercase\expandafter{\romannumeral2}}
\def\rmT{\uppercase\expandafter{\romannumeral3}}

\newtheorem*{conjecture}{Conjecture}

\def\const{\mathrm{const}}
\def\torus{\mathbb T}

\def\real{\mathbb R}
\def\integer{\mathbb Z}
\def\rational{\mathbb Q}
\def\bN{\mathbb N}

\newtheorem{theorem}{Theorem}
\newtheorem{prop}{Proposition}

\newtheorem{lemma}{Lemma}

\title{Geometry of plane sections\\ of the infinite regular skew polyhedron $\fsf$}

\author{Roberto De Leo, Ivan A. Dynnikov \thanks{INFN Cagliari,
Italy $\langle$roberto.deleo@ca.infn.it$\rangle$; Moscow State University, Russia
$\langle$dynnikov@mech.math.msu.su$\rangle$.
The work of the second named author is supported by Russian Federal Science Agency (grant
no.~{\cyr N\char'130}-1824.2008.1)}}

\begin{document}
\maketitle

\begin{abstract}
The asymptotic behavior of open plane sections of triply periodic surfaces is dictated,
for an open dense set of plane directions, by an integer second homology
class of the three-torus. The dependence of this homology class
on the direction can have a rather rich structure, leading in special cases to a fractal.
In this paper we present in detail the results for the skew polyhedron $\fsf$ and
in particular we show that in this case a fractal arises and that such a fractal
can be generated through an elementary algorithm, which in turn allows us to verify
for this case a conjecture of S.P.Novikov that such fractals have zero measure.
\end{abstract}

\pagestyle{myheadings}
\thispagestyle{plain}
\markboth{R. De Leo, I.A. Dynnikov}{Plane sections of the skew polyhedron $\fsf$}
\section{Introduction}

The study of plane sections of triply periodic surfaces in $\real^3$ was
initiated by S.P.Novikov in~\cite{Nov82} who raised a question whether
open (unbounded) components of such a section have some nice asymptotic
behavior. This was motivated by an application to conductivity theory.
A number of general theoretical results has been obtained since then
by A.Zorich~\cite{Zor84}, I.Dynnikov \cite{Dyn97,Dyn99},
and R.DeLeo~\cite{DeL03b,DeL05a}. For a number of surfaces, R.DeLeo
performed a numerical simulation, which confirmed the general
conclusions of the theory~\cite{DeL03a}.

In~\cite{DeL06} the main results were generalized to polyhedra.
Among the class of piecewise linear triply periodic closed
surfaces, the one of infinite
skew polyhedra~\cite{Cox37} is the most suitable for
numerical explorations of the geometry
of plane sections.
In~\cite{DeL06} the case of the regular skew polyhedron
$\cP$ of type $\sff$ was studied numerically in detail, showing that
the dependence of
open section's asymptotics on the plane direction keeps its rich
fractal-like structure also in the piecewise linear case.

In this work we present the results for the dual of $\cP$, namely the cubic polyhedron
$\cC=\fsf$~\cite{GSS02}; note that
both $\cP$ and $\cC$ are rough PL-approximations
of the smooth surface $\{\cos x^1+\cos x^2+\cos x^3=0\}$, which was itself studied numerically
in~\cite{DeL03a}.
It turns out that the $\cC$ case is rather noteworthy because its correspondent fractal
can be generated recursively through a simple algorithm, which, on one hand,
allowed us, for the first time,
to verify, in this concrete case, the conjecture~\cite{MN03} that
the set of exceptional directions has zero measure,
and, on the other hand, made
possible a systematic
comparison with the numerical data obtained through the NTC software library~\cite{NTC}.
%
\section{Topological structure of plane sections of triply periodic surfaces}
Let $M\subset\Tt$ be an embedded closed null homologous surface
in the torus $\Tt$, $H=(h_1,h_2,h_3)\in(\real^3)^*$
a covector. We denote by $\widehat M$ the preimage of $M$
under the projection $\real^3\rightarrow\Tt=\real^3/\integer^3$.
We consider the sections of $M$ by planes
$\langle H,x\rangle=\const$ (we call them \emph{$H$-sections})
and we are interested in the asymptotical behavior of
their unbounded regular connected components (if any).
Since only the direction of covector $H$ matters,
sometimes we shall treat $H$ as a point of the
projective plane $\real{\mathrm P}^2$.

In studying this question, the foliation $\mathcal F_H$
induced on $M$ by the closed one-form $\omega=(h_1dx^1+h_2dx^2+h_3dx^3)|_M$
plays the crucial role.
It is well known that, with probability $1$, in a proper sense,
a smooth closed one-form
whose critical points are all saddles induce dense leaves
on $\cM$. However, by restricting attention to a special
class of one-forms that are
pull-backs of a constant one-form on $\Tt$,
we fall exactly in the opposite situation, namely,
with probability $1$, open leaves are either absent or
confined to genus one minimal components of the foliation on $\cM$,
and dense leaves arise only in exceptional cases.

There are three principally different types of foliations
$\mathcal F_H$ and corresponding $H$-sections that may arise,
which we call \emph{trivial}, \emph{integrable}, and \emph{chaotic}.
Most typically, trivial means that all regular
leaves of $\mathcal F_H$ are closed,
integrable means that minimal components of $\mathcal F_H$ filled
with open leaves are of genus one, and in the chaotic case there is
a minimal component of genus $>1$. More precise definitions
are as follows.

Let $N\subset\torus^3$ be a piece-wise smooth embedded surface in $\torus^3$
such that $N\setminus M$ consists of disjoint open disks each of which lies
in a plane of the form $\{x\in\real^3\;;\;\langle H,x\rangle=\const\}$.
Such a surface is obtained by, first, cutting $M$ along some closed null homologous
leaves of $\mathcal F_H$ or null homologous saddle connection cycles,
second, removing some of the obtained connected components, and then gluing
up planar disks in order to obtain a closed surface.
For such a surface $N$, any leaf of $\mathcal F_H$ is either contained
in $N$ or disjoint from $N$. In the former case we say that the leaf
is \emph{absorbed} by $N$. When saying this we shall assume that $N$
is of the just specified form.

\begin{description}
\item[Trivial case.]
Every leaf of $\mathcal F_H$ is absorbed by a
two-sphere or a null homologous two-torus. If covector
$H$ is \emph{totally irrational}, i.e., $\dim_{\rational}\langle h_1,h_2,h_3\rangle
=3$, then this just means that all connected components of all $H$-sections
of $M$ are compact. If $\dim_{\rational}\langle h_1,h_2,h_3\rangle<3$,
then, additionally, \emph{periodic}, i.e., invariant under a non-trivial
shift, unbounded component of $H$-sections may arise;
\item[Integrable case.]
Every leaf of $\mathcal F_H$ is absorbed by a sphere or a two-torus,
and at least one leaf is absorbed
by a two-torus with non-zero homology class. In this case,
every regular non-closed component
of an $H$-section is a finitely deformed straight line,
i.e., it has the form $\gamma(t)=t\cdot v+O(1)$ for some parametrization,
where $v$ is a non-zero vector.
If $\dim_{\rational}\langle h_1,h_2,h_3\rangle=3$,
then the non-zero homology class of the tori absorbing
leaves of $\mathcal F_H$ is uniquely defined
up to sign. We denote it by $L_{M,H}$ and
consider as an integral covector in $\real^3$.
The identification of $H_2(\torus^3,\real)$ and $(\real^3)^*=
H^1(\torus^3,\real)$ is given by the Poincare duality.
This covector must obviously vanish at vector $v$: $\langle L_{M,H},v\rangle=0$.
So, if we assume our three-space Euclidean,
then (unless $H$ and $L_{M,H}$ are colinear)
we can simply write $v=L_{M,H}\times H$.
If $\dim_{\rational}\langle h_1,h_2,h_3\rangle<3$ the covector
$L_{M,H}$ may not be uniquely defined (up to sign),
but there may be at most two different choices.
We denote the projective class $(L_1:L_2:L_3)\in\rational\mathrm{P}^2$
of $L_{M,H}$ by $\ell_{M,H}$ and call \emph{the soul} of the foliation
$\mathcal F_H$.
\item[Chaotic case.] None of the above.
If $\dim_{\rational}\langle h_1,h_2,h_3\rangle=3$,
this means that a minimal component of $\mathcal F_H$
has genus $\geqslant3$.
The behavior of the corresponding $H$-sections have not been studied,
but the known examples suggest that, typically, a chaotic $H$-section
contains a single unbounded curve that ``wonders all around
the plane'', i.e., a $d$-neighborhood of the curve is the whole
plane for some finite $d$.
\end{description}

For a fixed surface $M$ and a rational point $\ell\in\rational\mathrm{P}^2
\subset\real\mathrm{P}^2$
we denote by $\zone_{M,\ell}$ the set
$$\zone_{M,\ell}=\{H\in\real\mathrm{P}^2\;;\;\ell_{M,H}=\ell\}.$$
If $\ell_{M,H}$ is not uniquely defined then the point $(h_1:h_2:h_3)$
is attributed to both corresponding subsets.
The set of points $(h_1:h_2:h_3)$ such that
the $H$-sections of $M$ are chaotic will be denoted by $\cE(M)$.

The following three propositions are extracted from~\cite{Dyn99}.

\begin{prop}\label{p1}
For a generic surface $M\subset\torus^3$
the sets $\zone_{M,\ell}$ are disjoint closed domains with piece-wise
smooth boundary. The set $\cE(M)$ is disjoint from $\rational\mathrm P^2$
and has zero measure. The set of
directions $H$ with trivial $H$-sections is open.
\end{prop}

In other words, the first claim of this proposition says that $\ell$, where defined,
is a locally constant function of $H$. We call the non-empty
domains $\zone_{M,\ell}$ \emph{stability zones} and refer
to $\ell$ as \emph{the label} of the stability zone $\zone_{M,\ell}$.

For studying the stability zones, it is usefull to consider
a 1-parametric family $M_c=\{x\in\torus^3\;;\;f(x)=c\}$ of
level surfaces of a fixed smooth function.

\begin{prop}\label{p2}
For a generic function $f$, there are continuous
functions $e_1,e_2:\real{\mathrm P}^2\rightarrow\real$ such that
\begin{itemize}
\item
$e_1(H)\leqslant e_2(H)$ for all $H\in\real{\mathrm P}^2$;
\item
the $H$-sections of $M_c$ are trivial if and only if $c\notin[e_1(H),e_2(H)]$;
\item
if $e_1(H)<e_2(H)$, then the $H$-sections of $M_c$ are integrable
for all $c\in[e_1(H),e_2(H)]$, and the soul $\ell$ of the corresponding
foliation $\mathcal F_{c,H}$ is independent of $c$.
\end{itemize}
\end{prop}

We define
\emph{generalized stability zones} $\zone_{f,\ell}$ as
$\mathcal D_{f,\ell}=\cup_c\mathcal D_{M_c,\ell}$,
and the set $\cE(f)$ as $\cE(f)=\cup_c\cE(M_c)$.

\begin{prop}\label{p3}
For a generic $f$,
generalized stability zones are closed
domains with piece-wise smooth boundary. If $\ell\ne\ell'$ then
the zones $\mathcal D_{f,\ell}$ ¨ $\mathcal D_{f,\ell'}$ can only have
intersections at the boundary, and, moreover, the number of their common
points is at most countable. If the whole $\real{\mathrm P}^2$
is not covered by a single generalized stability zone, then
the number of zones is countably infinite, and the set
$\cE(f)=\real{\mathrm P}^2\setminus(\cup_\ell\zone_{f,\ell})$
is non-empty and uncountable.
\end{prop}

It may happen that there is just one generalized stability zone
(say, a small enough perturbation of the function $\sin(x^1)$ will work),
but it is also easy
to find a function $f$ with non-empty $\cE(f)$: any function
with cubical symmetry is such.
In all examples known to us
two different generalized stability zones have at most one point in common.

It follows from Propositions~\ref{p1}--\ref{p3} that
the union $\cup_{\ell\in\rational\mathrm P^2}\text{int}(\mathcal D_{f,\ell})$
of the interiors of the zones is an open everywhere dense subset of
$\real\mathrm P^2$ and its complement $\overline{\mathcal E(f)}$
has the form of a two-dimensional cut out fractal set.

\begin{prop}[\cite{DeL03b}]\label{p4}
  If there is more than one generalized stability zone, then $\overline{\mathcal E(f)}$ is
  the set of accumulation points of the set of their souls.
\end{prop}

It is plausible but
still unknown whether $\mathcal E(f)$ has always zero measure.
The following stronger conjecture was proposed in~\cite{MN03}.

\begin{conjecture}
Whenever $\cE(f)$ is non-empty, the Hausdorff dimension of $\cE(f)$
is strictly between 1 and 2 for every $f$.
\end{conjecture}
So far only numerical checks of this conjecture were available
in the literature~\cite{DeL03a,DeL06};
in the next sections we shall provide, for the
particular case of the polyhedron $\cC$,
a full proof of the weaker, zero measure, conjecture and good numerical
evidence for the stronger one.
\section{Stability zones of $\cC$}
\label{sec:selfsim}
\begin{figure}
  \begin{center}
    \includegraphics[width=7.5cm]{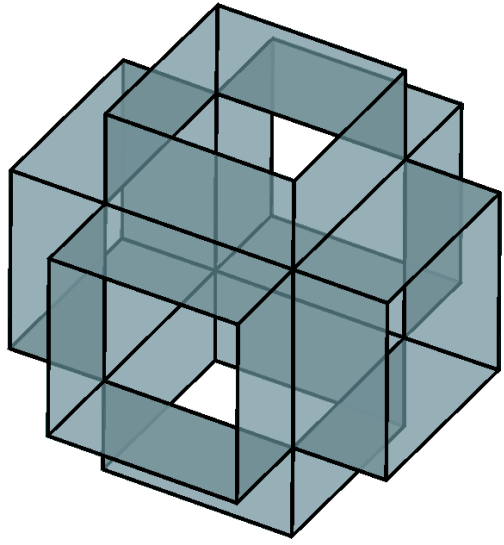}
  \end{center}
  \caption{%
    \small
    The $\fsf$ polyhedron embedded in the three-torus. Shown is the fundamental
    domain, which lies in the $[0,1]^3$ cube.
  }
  \label{fig:dc}
\end{figure}

The regular skew polyhedron $\cC=\fsf$ (see Figure~\ref{fig:dc}) is, up to isometries,
the unique cubic polyhedron with all monkey-saddle vertices~\cite{GSS02}; the vertices of its
embedding in $\Tt=[0,1]^3/\sim$ are the eight points in the orbit of $P=(1/4,1/4,1/4)$ under
the cubic symmetry group.
As level surface, $\cC$ can be represented in the $[0,1]^3$ cube as $\mf^{-1}(0)$ for
$$\mf(x^1,x^2,x^3)=\hbox{mid}(|2x^1-1|,|2x^2-1|,|2x^3-1|)-\frac12,$$
where $\hbox{mid}(a,b,c)$ is the middle value among $a$, $b$, and $c$.
\par
This surface, as well as the surface $\cos x^1+\cos x^2+\cos x^3=0$
and $\cC$'s dual---the truncated octahedron, has a very strong symmetry,
namely, its exterior is equal
to its interior modulo a translation. This means that for $\mf$ the functions $e_{1,2}$ mentioned
in Propositoin~\ref{p2} are such that $e_1=-e_2$, and hence,
stability zones of the surface $\cC$ coincide
with generalized stability zones of the function $\mf$.

Let us denote by $\psi_1,\psi_2,\psi_3$ the following projective transformations:
$$\begin{aligned}
\psi_1(h_1:h_2:h_3)&=(h_1:h_2+h_1:h_3+h_1),\\
\psi_2(h_1:h_2:h_3)&=(h_1+h_2:h_2:h_3+h_2),\\
\psi_3(h_1:h_2:h_3)&=(h_1+h_3:h_2+h_3:h_3).
\end{aligned}$$

\begin{theorem}\label{zones}
For the surface $\cC$ the stability zones are as follows:
$$\begin{aligned}
\mathcal D_{(1:0:0)}(\cC)&=\{(h_1:h_2:h_3)\in\real\mathrm P^2\;;\;h_1\geqslant
\left|h_2\right|+\left|h_3\right|\},\\
\mathcal D_{(1:1:1)}(\cC)&=\{(h_1:h_2:h_3)\in\real\mathrm P^2\;;\;
\left|h_1\right|+\left|h_2\right|+\left|h_3\right|\leqslant 4h_1,4h_2,4h_3\},\\
\mathcal D_{\psi_{i_1}(\psi_{i_1}(\ldots\psi_{i_k}((1:1:1))\ldots))}(\cC)&=
\psi_{i_1}(\psi_{i_1}(\ldots\psi_{i_k}(\mathcal D_{(1:1:1)}(\cC))\ldots)),
\end{aligned}$$
where $(i_1,\ldots,i_k)$ is an arbitrary finite sequence
of elements from $\{1,2,3\}$, and, in addition, all zones
obtained from the listed ones by cubical symmetries:
permutations and changing signs of the coordinates.
\end{theorem}

The proof of this theorem will rest on the following two lemmas.

\begin{lemma}\label{IDL1}
We have
$$\mathcal D_{(1:0:0)}(\cC)\supset\{(h_1:h_2:h_3)\in\real\mathrm P^2\;;\;h_1\geqslant
\left|h_2\right|+\left|h_3\right|\}.$$
\end{lemma}

\begin{figure}
\centerline{\includegraphics[height=7cm]{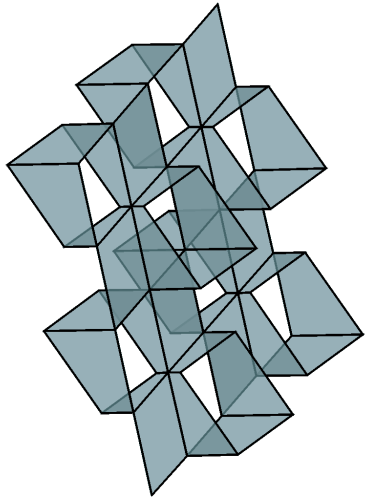}}
\caption{A connected component of $\widehat\cC$ after cutting
along closed leaves}\label{fig:wrapped}
\end{figure}

\begin{proof}
If the inequality
$h_1\geqslant
\left|h_2\right|+\left|h_3\right|$ is satisfied then
the plane $h_1x^1+h_2x^2+h_3x^3=\mathrm{const}$
passing through the center of the unit cube $[0,1]^3$
separates the faces $x^1=0$ and $x^1=1$.
This means that 
closed leaves of the corresponding
foliation will cut our triply periodic surface $\widehat\cC$
into parts each of which is a finitely deformed
plane $x^1=\mathrm{const}$ with holes (see Figure~\ref{fig:wrapped}).
Filling the holes by flat disks and projecting to $\torus^3$ we
obtain two tori whose homology class is equal up to sign to
$(1,0,0)\in(\real^3)^*$.
\end{proof}

\begin{lemma}\label{IDL2}
Let $h_1,h_2,h_3\geqslant0$, $(h_1,h_2,h_3)\ne(0,0,0)$. Then we have
$(h_1:h_2:h_3)\in\mathcal D_\ell$ if and only if
$\psi_i((h_1:h_2:h_3))\in\mathcal D_{\psi_i(\ell)}$,
where $i=1,2,3$.
\end{lemma}

\begin{proof}
For convenience, we shift the coordinate system as follows:
$(x^1,x^2,x^3)\mapsto(x^1+1/4,x^2+1/4,x^3+1/4)$. Our surface
$\widehat\cC$ cuts $\real^3$ into two parts, $N_-$ and $N_+$, that
can now be characterized as follows: $N_-$ (resp.\ $N_+$) consists of
points $(x^1,x^2,x^3)\in\real^3$ such that
at least two (resp.\ at most one) of the three numbers $\{x^1\},\{x^2\},\{x^3\}$
are in the interval $[0,1/2]$ (resp.\ $(0,1/2)$), where $\{x\}$ denotes
the fractional part of $x$.

We may assume $i=3$, $(h_1,h_2,h_3)=(\alpha,\beta,1)$ without loss of generality.
Let $\Pi$ be the plane defined by $\alpha x^1+\beta x^2+x^3=c$.
Denote by $Q_-$ the projection of $\Pi\cap N_-$ to the $x^1,x^2$-plane
along $x^3$. According to the description of $N_-$ given above
$Q_-$ is the set of points $(x^1,x^2)\in\real^2$
such that exactly two of the three numbers $\{x^1\},\{x^2\},\{c-\alpha x^1-\beta x^2\}$
are in the interval $[0,1/2]$.

Denote by $\square_{a,b}$ the square $\{(x^1,x^2)\in\real^2\,;\,
a\leqslant x^1\leqslant a+1/2,\,
b\leqslant x^2\leqslant b+1/2\}$, and by $S_m$ the strip
defined by $m\leqslant c-\alpha x^1-\beta x^2\leqslant m+1/2$.
We have
$$Q_-=\Bigl(\bigcup_{j,k\in\integer}\square_{j,k}\Bigr)
\cup\Bigl(\bigcup_{j,k,m\in\integer}\square_{j+1/2,k}\cap S_m\Bigr)
\cup\Bigl(\bigcup_{j,k,m\in\integer}\square_{j,k+1/2}\cap S_m\Bigr).$$
The first part in this union, $\bigcup_{j,k\in\integer}\square_{j,k}$,
does not depend on $\Pi$. We call these squares \emph{mainlands}.

Each intersection $\square_{j+1/2,k}\cap S_m$ with $j,k,m\in\integer$,
whenever non-empty, is a convex polygon that can be of the following three types:
\begin{description}
\item[\emph{cape}:] it has a piece of boundary in common with exactly
one of the mainlands
$\square_{j,k}$ or $\square_{j+1,k}$;
\item[\emph{bridge}:] it has a piece of boundary in common with both mainlands
$\square_{j,k}$ and $\square_{j+1,k}$;
\item[\emph{island}:] it is disjoint from mainlands.
\end{description}
Similarly, one defines the type of a polygon $\square_{j,k+1/2}\cap S_m$
regarding adjacency to the mainlands $\square_{j,k}$ and $\square_{j,k+1}$,
see Figure~\ref{fig:geography}.

\begin{figure}
\centerline{%
\begin{picture}(240,230)
\put(0,30){\includegraphics[height=200pt]{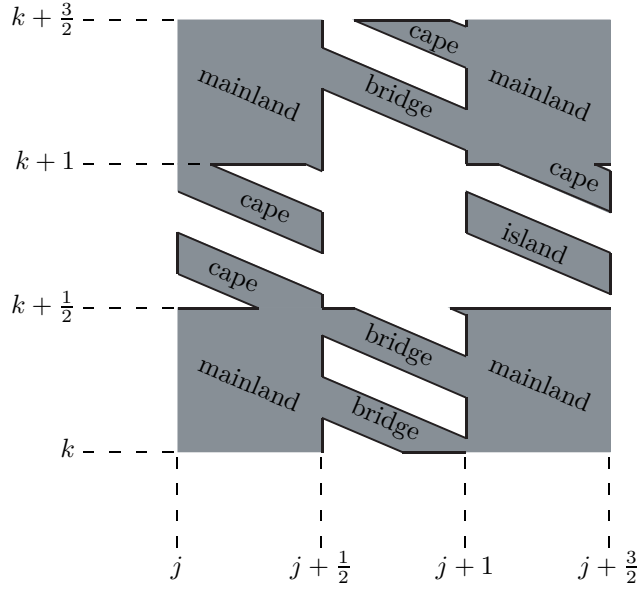}}
\put(130,64){\rotatebox{-23}{bridge}}
\put(135,92){\rotatebox{-23}{bridge}}
\put(135,186){\rotatebox{-23}{bridge}}
\put(73,80){\rotatebox{-23}{mainland}}
\put(183,80){\rotatebox{-23}{mainland}}
\put(183,190){\rotatebox{-23}{mainland}}
\put(73,190){\rotatebox{-23}{mainland}}
\put(76,117){\rotatebox{-23}{cape}}
\put(88,143){\rotatebox{-23}{cape}}
\put(152,207){\rotatebox{-23}{cape}}
\put(205,155){\rotatebox{-23}{cape}}
\put(188,130){\rotatebox{-23}{island}}
\multiput(30,49)(10,0)4{\line(1,0)4}
\put(27,47){\hbox to0pt{\hss$k$}}
\multiput(30,103.5)(10,0)4{\line(1,0)4}
\put(27,101.5){\hbox to0pt{\hss$k+\frac12$}}
\multiput(30,158)(10,0)5{\line(1,0)4}
\put(27,156){\hbox to0pt{\hss$k+1$}}
\multiput(30,212.5)(10,0)4{\line(1,0)4}
\put(27,210.5){\hbox to0pt{\hss$k+\frac32$}}
\multiput(65.5,13)(0,10)4{\line(0,1)4}
\put(65.5,3){\hbox to 0pt{\hss$j$\hss}}
\multiput(120,13)(0,10)4{\line(0,1)4}
\put(120,3){\hbox to 0pt{\hss$j+\frac12$\hss}}
\multiput(174.5,13)(0,10)4{\line(0,1)4}
\put(174.5,3){\hbox to 0pt{\hss$j+1$\hss}}
\multiput(229,13)(0,10)4{\line(0,1)4}
\put(229,3){\hbox to 0pt{\hss$j+\frac32$\hss}}
\end{picture}%
}
\caption{Geography of the section $\Pi\cap N_-$}\label{fig:geography}
\end{figure}

Capes are not interesting for us because their removal is equivalent to
a finite deformation of $Q_-$. It is not hard to write
the necessary and sufficient condition for $\square_{j+1/2,k}\cap S_m$
to be a bridge:
\begin{equation}\label{bridge}
c-\alpha\Bigl(j+\frac12\Bigr)-\beta\Bigl(k+\frac12\Bigr)-\frac12
\leqslant m\leqslant c-\alpha(j+1)-\beta k,
\end{equation}
or an island:
\begin{equation}\label{island}
c-\alpha\Bigl(j+\frac12\Bigr)-\beta\Bigl(k+\frac12\Bigr)-\frac12
\geqslant m\geqslant c-\alpha(j+1)-\beta k.
\end{equation}

Now we apply $\psi_3$ to $H$, which gives $H'=(\alpha',\beta',1)=(\alpha+1,\beta+1,1)$.
Let $Q_-'$ and $S_m'$ be defined in the same way as $Q_-$ and $S_m$
with $\alpha$ and $\beta$
replaced by $\alpha'=\alpha+1$ and $\beta'=\beta+1$, respectively.
Then the intersection $\square_{j+1/2,k}\cap S_m'$ is a bridge if and only if
$$c-(\alpha+1)\Bigl(j+\frac12\Bigr)-(\beta+1)\Bigl(k+\frac12\Bigr)-\frac12
\leqslant m\leqslant c-(\alpha+1)(j+1)-(\beta+1)k,$$
which can be rewritten as
$$c-\alpha\Bigl(j+\frac12\Bigr)-\beta\Bigl(k+\frac12\Bigr)-\frac12
\leqslant m+j+k+1\leqslant c-\alpha(j+1)-\beta k.$$
Thus, $\square_{j+1/2,k}\cap S_m$ is a bridge if
and only if so is $\square_{j+1/2,k}\cap S_{m-j-k-1}'$.
Similarly, the same
is true about islands as well as bridges and islands in
squares of the form $\square_{j,k+1/2}$.

So, bridges and islands of $Q_-'$ in the square $\square_{j+1/2,k}$
or $\square_{j,k+1/2}$ are in a natural one to one correspondence with
those of $Q_-$.
Therefore, $Q_-'$ and $Q_-$ are obtained from each other by
a finite deformation. The geometrical difference between
$Q_-'$ and $Q_-$ can be vaguely described as follows: islands
and bridges of $Q_-'$ are narrower than those of $Q_-$,
and $Q_-'$ has more capes. See Figure~\ref{fig:morecapes} for an example.

\begin{figure}
$$\begin{array}{cc}
\includegraphics[height=4cm]{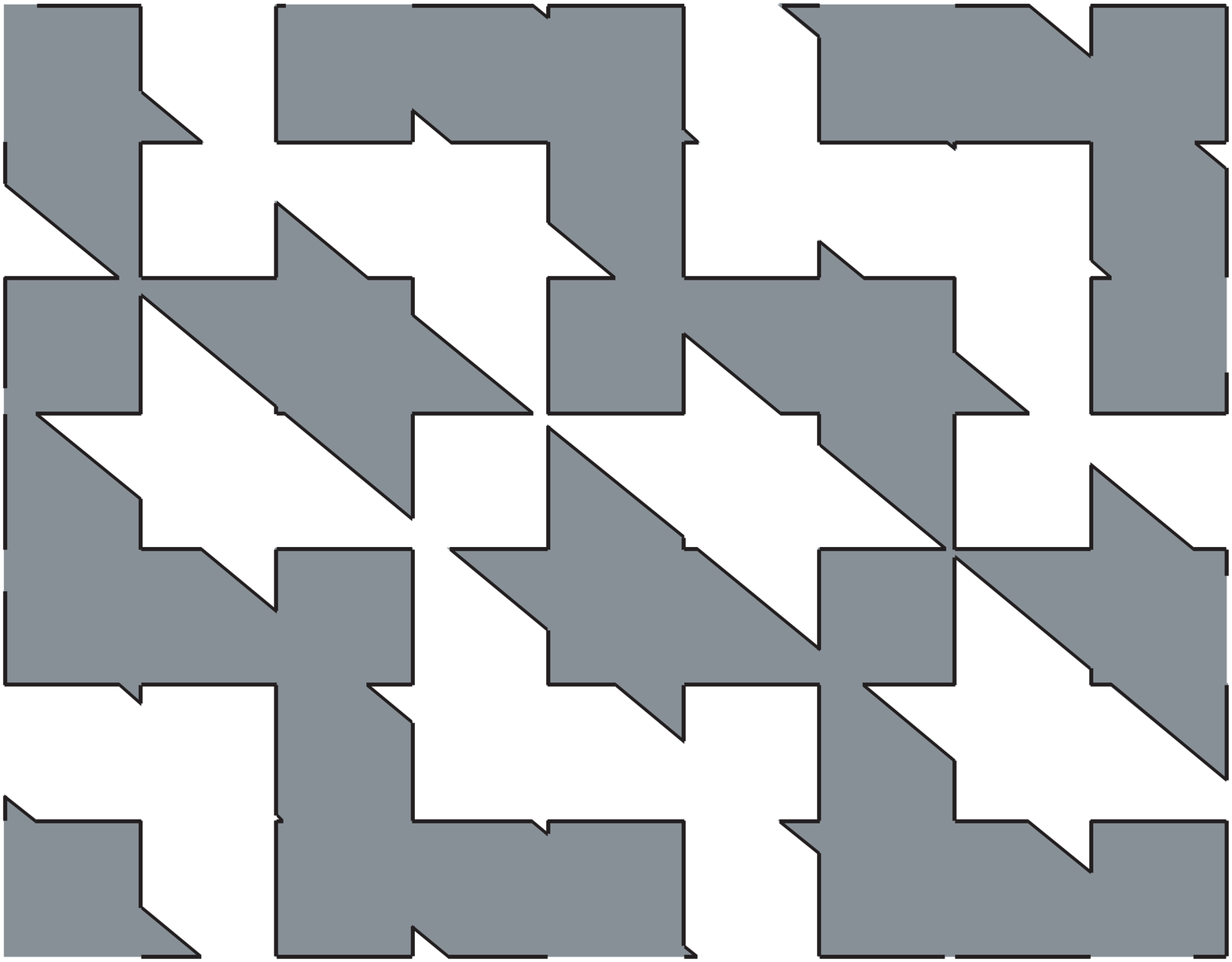}&
\includegraphics[height=4cm]{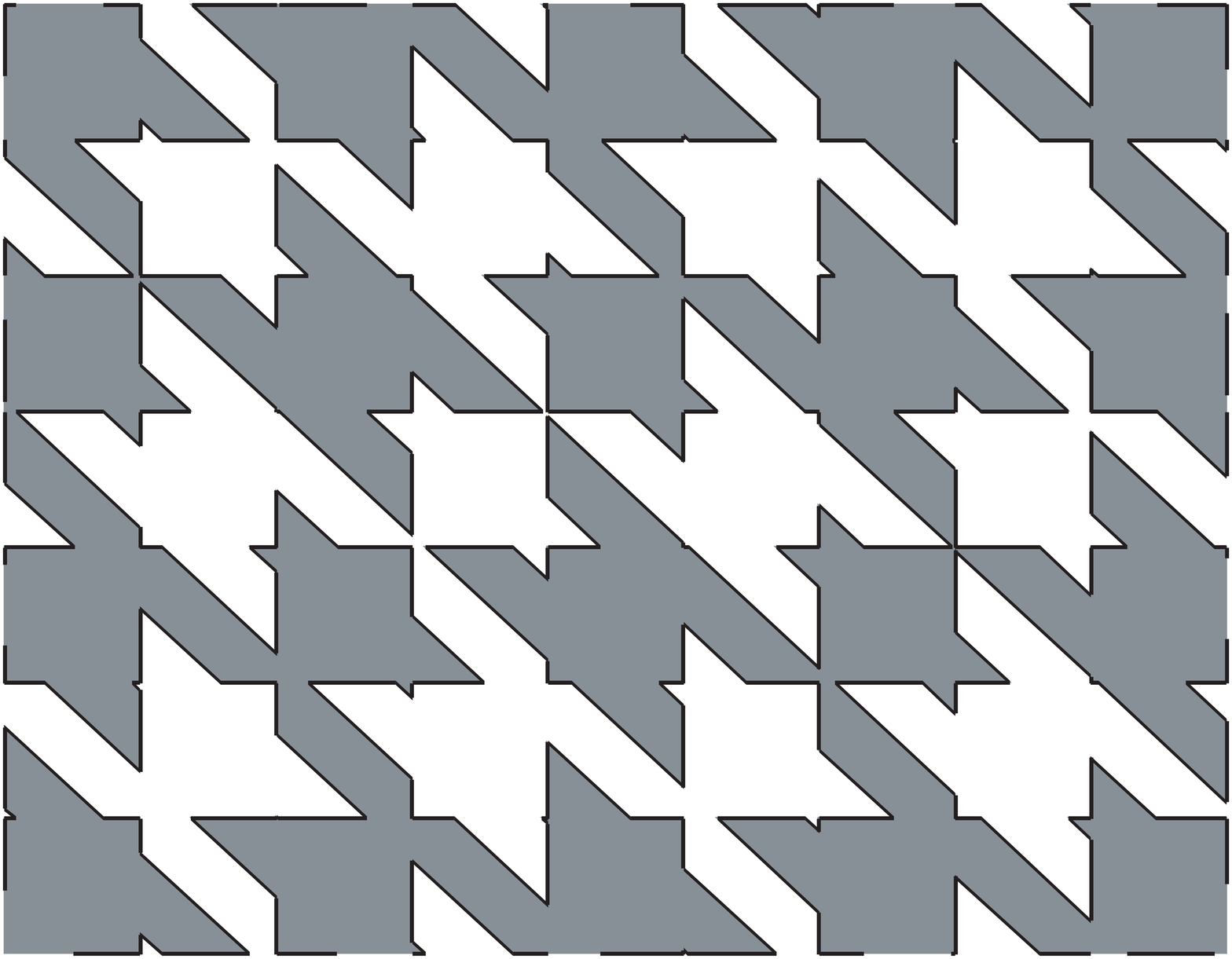}\\
Q_-&Q_-'
\end{array}$$
\caption{The transition from $Q_-$ to $Q_-'$ produces more capes and makes
the bridges narrower}\label{fig:morecapes}
\end{figure}

In the genus three case the integrability of our foliation
is equivalent to the existence of closed fibres of the foliation (or null homologous
saddle connection cycles).
We have just seen that $H$-sections and $H'$-sections are obtained
from each other by a finite deformation. Hence, they are both integrable
or both chaotic. If they are integrable, let $\ell$ and $\ell'$
be the labels of the corresponding zones. We want to show
that $\ell'=\psi_3(\ell)$. Since both $\ell$ and $\ell'$
are locally constant functions of $H$ it is enough to consider
the case of totally irrational $H$. Then the asymptotic direction
$v=H\times\ell$ is of irrationality degree two, i.e., $\dim_\rational\langle
v^1,v^2,v^3\rangle=2$, and $\ell$
is the only rational covector up to multiple that
vanishes at $v$. Thus, it suffices to show
that $\psi_3(\ell)$ vanishes at $v'$, the direction
of open components of $H'$-sections.

The latter follows easily from the fact that the projections of $v$ and
$v'$ to the $x^1,x^2$-plane coincide, which, in turn, follows from
the coincidence, up to finite deformation, of $Q_-$ and $Q_-'$.
\end{proof}

\begin{proof}[Proof of Theorem~\ref{zones}]
Due to the cubical symmetry of the surface it suffices
to establish the fractal structure in the region
$C_+=\{(h_1:h_2:h_3)\in\real\mathrm P^2\,;\,h_1,h_2,h_3\geqslant0\}$.
Put
$$\begin{aligned}
\mathcal R_{(1:0:0)}&=\{(h_1:h_2:h_3)\in\real\mathrm P^2\;;\;h_1\geqslant
\left|h_2\right|+\left|h_3\right|\},\\
\mathcal R_{(0:1:0)}&=\{(h_1:h_2:h_3)\in\real\mathrm P^2\;;\;h_2\geqslant
\left|h_1\right|+\left|h_3\right|\},\\
\mathcal R_{(0:0:1)}&=\{(h_1:h_2:h_3)\in\real\mathrm P^2\;;\;h_3\geqslant
\left|h_1\right|+\left|h_2\right|\},\\
\mathcal R_{(1:1:1)}&=\{(h_1:h_2:h_3)\in\real\mathrm P^2\;;\;
\left|h_1\right|+\left|h_2\right|+\left|h_3\right|\leqslant 4h_1,4h_2,4h_3\},\\
\mathcal R_{\psi_{i_1}(\psi_{i_1}(\ldots\psi_{i_k}((1:1:1))\ldots))}&=
\psi_{i_1}(\psi_{i_1}(\ldots\psi_{i_k}(\mathcal R_{(1:1:1)})\ldots)),
\end{aligned}$$
and $\mathcal R_\ell=\varnothing$ if $\ell\notin\{(1:0:0),(0:1:0),(0:0:1)\}$
and $\ell$ is not of the form
$\psi_{i_1}(\psi_{i_1}(\ldots\psi_{i_k}((1:1:1))\ldots))$.
We want to show that $\mathcal R_\ell=\mathcal D_\ell$ for all $\ell\in
\rational\mathrm P^2\cap C_+$. We have already shown in Lemma~\ref{IDL1}
that $\mathcal R_\ell\subset\mathcal D_\ell$ for $\ell=(1:0:0)$
By symmetry it is also true for $\ell=(0:1:0),(0:0:1)$.

It is a straightforward check to see that
$$\mathcal R_{(1:1:1)}=\psi_1(\mathcal R_{(1:0:0)}\cap C_+)\cup
\psi_2(\mathcal R_{(0:1:0)}\cap C_+)\cup\psi_3(\mathcal R_{(0:0:1)}\cap C_+).$$
By Lemma~\ref{IDL2} this implies
$\mathcal R_{(1:1:1)}\subset\mathcal D_{(1:1:1)}$
as $(1:1:1)=\psi_1(1:0:0)=\psi_2(0:1:0)=\psi_3(0:0:1)$,
and, therefore,
$\mathcal R_\ell\subset\mathcal D_\ell$ for all $\ell\in C_+\cap\rational\mathrm P^2$.

In order to establish Theorem~\ref{zones} it suffices to show that
the zones $\mathcal D_\ell$ are not larger that $\mathcal R_\ell$, and
there are no other stability zones. Both claims follow from
the fact that $\bigcup_\ell\mathcal R_\ell$ covers all rational points:
$$\bigcup_\ell\mathcal R_\ell\supset C_+\cap\rational\mathrm P^2.$$

Indeed, let $\varphi$ be the following map from $C_+$ to itself:
$$\varphi(h_1:h_2:h_3)=\left\{\begin{aligned}
(h_1:h_2-h_1:h_3-h_1)&,\text{ if }0\leqslant h_1\leqslant h_2,h_3,\\
(h_2:h_3:h_1)&,\text{ otherwise}.
\end{aligned}
\right.$$
By construction, for every $H\in C_+$ we have
$H\in\bigcup_\ell\mathcal R_\ell$ if and only if
$\varphi(H)\in\bigcup_\ell\mathcal R_\ell$.
If $H$ is a rational covector from $C_+$, then after applying
$\varphi$ finitely many times, one of the coordinates of
the obtained covector becomes zero. All such points are
covered by $\mathcal R_{(1:0:0)}$, $\mathcal R_{(0:1:0)}$
and $\mathcal R_{(0:0:1)}$.
\end{proof}
\begin{figure}
  \begin{center}
    \includegraphics[width=7cm]{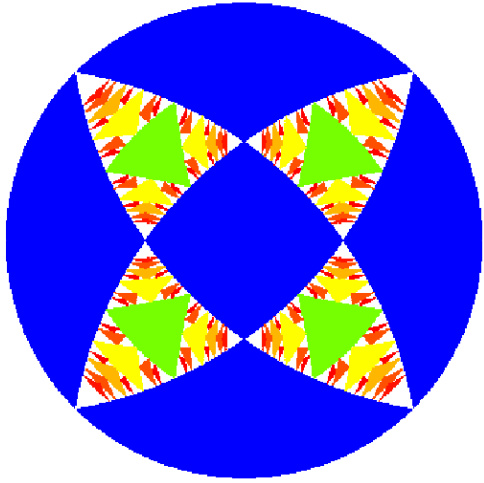}
    \caption{%
      \small
      Picture of the fractal in the disc model of $\RPt$. The center of the disc corresponds to the $z$ axis,
      so the central square is the stability zone $\cD_{(0:0:1)}$, the one touching it in the left and right vertices is
      $\cD_{(1:0:0)}$ and the third one is $\cD_{(0:1:0)}$. The green triangles are the
      stability zones $\cD_{(1:1:1)}$,
      $\cD_{(1:1:-1)}$, $\cD_{(1:-1:1)}$ and $\cD_{(1:-1:-1)}$.
    }
    \label{fig:frDisc}
  \end{center}
\end{figure}
So, we have the following picture in $\real\mathrm P^2$: four lines
$h_1\pm h_2\pm h_3=0$ cut $\real\mathrm P^2$ into
three ``squares'', which are zones $\mathcal D_{(1:0:0)}$,
$\mathcal D_{(0:1:0)}$, $\mathcal D_{(0:0:1)}$, and
four triangles obtained from each other by cubical symmetries,
in which there are infinitely many stability zones.
We shall concentrate on the triangle that is contained in $C_+$.
This triangle has vertices $(1:1:0)$, $(1:0:1)$, $(0:1:1)$
and is defined by the inequalities
$$\frac{h_1+h_2+h_3}2\geqslant h_1,h_2,h_3\geqslant0.$$
Let us denote this triangle by $\Delta$.
The zone $\mathcal D_{(1:1:1)}$ is also a triangle that is contained
in $\Delta$ and has its vertices, $(2:1:1)$, $(1:2:1)$, $(1:1:2)$,
at the sides of $\Delta$.
The complement $\Delta\setminus\mathcal D_{(1:1:1)}$
consists of three triangles that are exactly $\psi_1(\Delta)$,
$\psi_2(\Delta)$, and $\psi_3(\Delta)$.
In each triangle $\Delta_1=\psi_1(\Delta)$, $\Delta_2=\psi_2(\Delta)$,
$\Delta_3=\psi_3(\Delta)$ the picture is obtained from
that in $\Delta$ by the corresponding projective transformation $\psi_i$.

For a finite sequence $a=(a_1,a_2,\dots,a_k)$ of indices $1,2,3$
we denote by $\psi_a$ the mapping $\psi_{a_1}\circ\psi_{a_2}\circ\ldots
\circ\psi_{a_k}$. By $a'$, $a''$, and $a'''$ we denote
the sequences $(a_1,\dots,a_k,1)$, $(a_1,\dots,a_k,2)$, $(a_1,\dots,a_k,3)$,
respectively. For any such sequence $a$ we have the following.
The triangle $\Delta_a=\psi_a(\Delta)$
is bounded by the zones $\mathcal D_{\psi_a(1:0:0)}$,
$\mathcal D_{\psi_a(0:1:0)}$, $\mathcal D_{\psi_a(0:0:1)}$.
It contains the zone
$\mathcal D_{\psi_a(1:1:1)}$ whose vertices
$\psi_a(2:1:1)$, $\psi_a(1:2:1)$, $\psi_a(1:1:2)$
are on the sides of $\Delta_a$, and the complement
$\Delta_a\setminus\mathcal D_{\psi_a(1:1:1)}$ is the union
of the triangles $\Delta_{a'}$, $\Delta_{a''}$, $\Delta_{a'''}$.

\begin{figure}
  \begin{center}
    \includegraphics[width=7.5cm,viewport=30 130 580 668,clip]{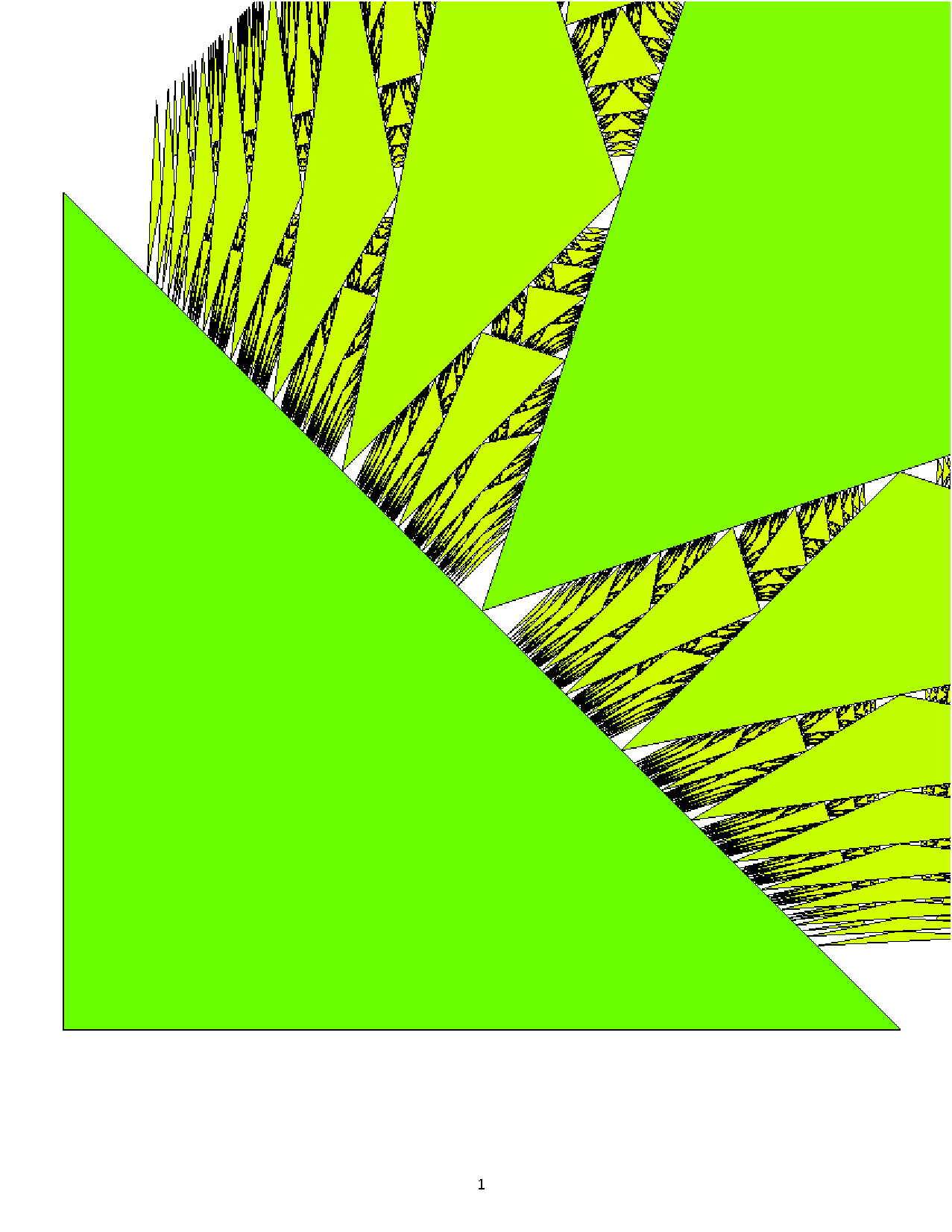}
    \caption{%
      \small
      Picture of the fractal in the square $[0,1]^2$ in the $z=1$ projective chart.
    }
    \label{fig:frDiscAn}
  \end{center}
\end{figure}

\begin{prop}
The intersection $\mathcal E(\cC)\cap C_+$
consists of points of the form
$$\lim_{k\rightarrow\infty}\psi_{i_1}(\psi_{i_2}(\ldots
\psi_{i_k}((1:1:1))\ldots)),$$
where $(i_1,i_2,\ldots)$ runs over all possible sequences
of elements from $\{1,2,3\}$ containing each index infinitely many times.
Other points in $\mathcal E(\cC)$ are obtained from these by cubical simmetries.
\end{prop}

\begin{proof}
From the structure of stability zones established above it
follows that the intersection $\mathcal E\cap C_+$
is the union of subsets
$$\bigcap\limits_k\psi_{i_1}(\psi_{i_2}(\ldots
\psi_{i_k}(\Delta)\ldots))$$
over all sequences $(i_k)\in3^{\mathbb N}$ in which all three indices appear
infinitely many times.
It suffices to show that every such a subset is actually a single point,
which follows from Proposition~\ref{zoneasympt} below.
\end{proof}
\subsection{Measure of $\cE$}
%
In all cases studied so far no algorithm was found to generate all stability
zones and nothing could be
said about the measure of the $\cE(f)$.
This is therefore the first case in which it is possible to
check the truth of the zero measure conjecture.
\begin{theorem}
The set $\mathcal E(\cC)$ of ``chaotic'' directions has measure zero.
\end{theorem}
\begin{proof}
Again, due to the symmetry, it suffices to prove the claim for $\mathcal E\cap C_+$.
We denote this set by $\mathcal E_+$. As we have seen, it is contained in
the triangle $\Delta$, and the following holds:
$$\mathcal E_+=\psi_1(\mathcal E_+)\cup\psi_2(\mathcal E_+)\cup
\psi_3(\mathcal E_+).$$
Applying this recursively, we get
\begin{multline}
\mathcal E_+=\psi_1(\mathcal E_+)\cup\psi_2(\mathcal E_+)\cup
\psi_{31}(\mathcal E_+)\cup\psi_{32}(\mathcal E_+)\cup
\psi_{331}(\mathcal E_+)\cup\psi_{332}(\mathcal E_+)\cup\\
\ldots\cup
\psi_{\underbrace{\scriptstyle3\ldots3}_k1}(\mathcal E_+)\cup
\psi_{\underbrace{\scriptstyle3\ldots3}_k2}(\mathcal E_+)\cup
\psi_{\underbrace{\scriptstyle3\ldots3}_{k+1}}(\mathcal E_+).
\end{multline}

Let $\mu$ be the measure of $\mathcal E_+$ and $\mu_k$
the measure of $\psi_{\underbrace{\scriptstyle3\ldots3}_k1}(\mathcal E_+)$,
which is equal to the measure of
$\psi_{\underbrace{\scriptstyle3\ldots3}_k2}(\mathcal E_+)$.
The measure of the triangle
$\psi_{\underbrace{\scriptstyle3\ldots3}_{k+1}}(\mathcal E_+)$
tends to zero as $k$ goes to $\infty$. Therefore, we have
$$\mu=2\sum_{k=0}^\infty\mu_k.$$
The idea now is to show that $\mu_k\leqslant c_k\mu$, where $c_0,c_1,\ldots$
are constants such that
$$\sum_{k=0}^\infty c_k<\frac12,$$
which, together with the previous equality, implies $\mu=0$.

Now we provide the necessary technical details. First of all,
we need to parametrise the triangle $\Delta$.
We use the following parametrization: $(u,v)\mapsto
(1-v:1-u:u+v)$, $u,v\geqslant0$, $u+v\leqslant1$.

The property of a set to be of zero measure does not
depend on the choice of a regular Borel measure.
In order to get the desired inequalities we use a somewhat artificial
measure on $\Delta$, namely the following one:
$$\frac{du\,dv}{(1+u+v)^2}.$$
By doing so we keep the symmetry between $\psi_1$ and $\psi_2$ (one is conjugate
to the other by the permutation $u\leftrightarrow v$), so the measures
of the sets $\psi_{\underbrace{\scriptstyle3\ldots3}_k1}(\mathcal E_+)$
and $\psi_{\underbrace{\scriptstyle3\ldots3}_k2}(\mathcal E_+)$ coincide.
We denote by $f_k$ the mapping $\psi_{\underbrace{\scriptstyle3\ldots3}_k1}\circ R$,
where $R(h_1:h_2:h_3)=(h_3:h_1:h_2)$, written in the $u,v$-parametrization:
$$f_k(u,v)=\Bigl(\frac v{u+(k+1)(v+1)},\frac 1{u+(k+1)(v+1)}\Bigr).$$
Since we have $R(\mathcal E)=\mathcal E$, the image $f_k(\mathcal E_+)$
coincides with $\psi_{\underbrace{\scriptstyle3\ldots3}_k1}(\mathcal E_+)$.

Obviously, the measure of $f_k(\mathcal E_+)$ is bounded from above
by $c_k\mu$, where
$$c_k=\max_{u,v\geqslant0,\,u+v\leqslant1}
J_{f_k}(u,v)\frac{(1+u+v)^2}
{\Bigl(1+\frac v{u+(k+1)(v+1)}+\frac 1{u+(k+1)(v+1)}\Bigr)^2}.$$
By $J_f$ we denote the Jacobian of the mapping $f$. A routine check gives:
\begin{multline*}
c_k=\max_{u,v\geqslant0,\,u+v\leqslant1}
\frac{(1+u+v)^2}
{(u+(k+2)(v+1))^2(u+(k+1)(v+1))}=\\
\left\{\begin{array}{ll}
\displaystyle\frac14&\text{if }k=0,\\
\displaystyle\frac4{(k+2)(k+3)^2}&\text{otherwise},
\end{array}\right.
\end{multline*}
$$\sum_{k=0}^\infty c_k=\frac14+\sum_{k=1}^\infty\frac4{(k+2)(k+3)^2}=
\frac{253}{36}-\frac23\pi^2\approx0.448<\frac12,
$$
which completes the proof.

\end{proof}
\subsection{Asymptotics of stability zones and fractal dimension of $\cE$}
\label{sec:asympt}
\begin{figure}
  \begin{center}
    \includegraphics[width=7cm]{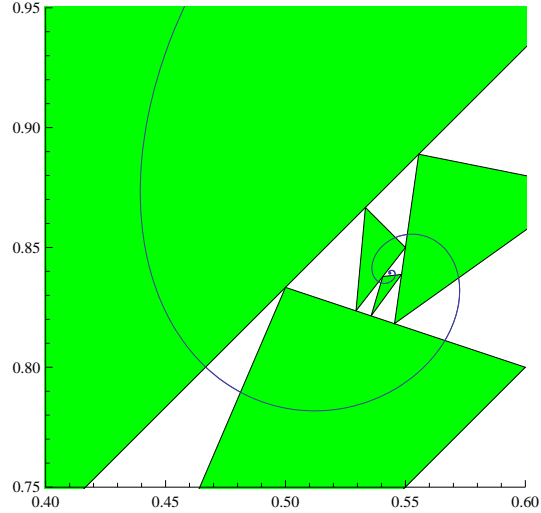}
    \caption{%
      \small
      Detail, in the $z=1$ projective chart, of the first few zones of the Tribonacci
      sequence of stability zones starting by $\cD_{(1:0:0)}$, $\cD_{(0:1:0)}$ and $\cD_{(0:0:1)}$.
      The labels of the zones shown above are, in decreasing area order, $(1:2:2)$, $(2:3:4)$,
      $(4:6:7)$, $(7:11:13)$ and $(13:20:24)$. All centers of the zones of the sequence lie on the
      smooth ``Tribonacci projective spiral'' drawn above which is converging to
      $(\alpha^2-\alpha-1,\alpha-1)\simeq(.543,.839)$.
    }
    \label{fig:spiral}
  \end{center}
\end{figure}

From what said at the beginning of the section it follows immediately that
all stability zones, except the biggest ones, which are squares, are triangles
and their labels satisfy a simple recursive relation.
\begin{prop}
\label{prop:alg}
  Every stability zone of $\cC$, except for the three squares with souls
  $(1:0:0)$, $(0:1:0)$, and $(0:0:1)$,
  are triangles; these triangles are generated starting from $\Delta$ (and its three symmetric
  triangles with respect to the coordinate planes) according to the the following recursive algorithm:
  \item{1.} evaluate the vector sums $w_1=v_2+v_3$, $w_2=v_3+v_1$, $w_3=v_1+v_2$ of all pairs
    of the three vertices $v_i$ of $\Delta$;
  \item{2.} consider the $w_i$ as projective coordinates and add the triangle having those points
    as vertices to the list of all stability zones -- its soul is given by the vector sum of the three vertices
    of $\Delta$;
  \item{3.} consider now the three triangles $\Delta_1=\langle v_1,w_2,w_3\rangle$,
  $\Delta_2=\langle w_1,v_2,w_3\rangle$,
    $\Delta_3=\langle w_1,w_2,v_3\rangle$ and repeat the algorithm for each of them.
\end{prop}
This fact can be exploited to find an explicit expression for elements of $\cE$ by considering
``spiralling down'' (or ``steepest descent'') sequences of stability zones (see Figure~\ref{fig:spiral}).
Select indeed an ordered triple
$(\alpha,\beta,\gamma)$ of stability zones which bind a triangle, namely such that two touch each other
and the third touches both, and build out of them the
recursive sequence of stability zones such that every
new stability zone is the one whose vertices touch the previous three stability zones. It is easy to see that the sequence
of these stability zones are arranged in a sort of spiral and by construction, the label
of every stability zone of this sequence is the sum of the labels of the previous three zones, so that
the sequence of the labels is a Tribonacci sequence in $\mathrm PH_2(\Tt,\bZ)$.
\begin{prop}
  The limit of every such sequences belongs to $\cE$.
  In particular $(\alpha^2-\alpha-1:\alpha-1:1)\in\cE$, where $\alpha$ is the Tribonacci
  constant, namely the real solution of the Tribonacci equation $x^3-x^2-x-1=0$.
\end{prop}
\begin{proof}
  All vertices of stability zones of $\cC$ are 1-rational points of $\RPt$ and therefore points on the
  boundaries have irrationality degree not higher than 2. Now, consider the sequence starting
  from the three squares $a_1=\cD_{(1:0:0)}$, $a_2=\cD_{(0:1:0)}$, $a_3=\cD_{(0:0:1)}$, so that the
  first few next terms of the sequence will be $a_4=\cD_{(1:1:1)}$, $a_5=\cD_{(1:2:2)}$, $a_6=\cD_{(2:3:4)}$
  and so on: a simple calculation show that the label of the $n$th stability zone of the sequence,
  modulo terms in $\beta^n$ and $\bar\beta^n$, will be
  $$
  (\frac{\beta\bar\beta}{(\alpha-\beta)(\alpha-\bar\beta)}\alpha^n
  :
  \frac{-\beta-\bar\beta}{(\alpha-\beta)(\alpha-\bar\beta)}\alpha^n
  :
  \frac{1}{(\alpha-\beta)(\alpha-\bar\beta)}\alpha^n)
  $$
  where $\alpha=(1+\sqrt[3]{19-3\sqrt{33}}+\sqrt[3]{19+3\sqrt{33}})/3\simeq1.84$ is the Tribonacci
  constant and $\beta,\bar\beta$ the remaining two solutions of the Tribonacci equation $x^3-x^2-x-1=0$.
  Since $|\beta|<\alpha$ the sequence of labels converges to the triple of coefficients of
  $\alpha^n$, namely $\ell_\infty=(\beta\bar\beta:-\beta-\bar\beta:1)=(\alpha^2-\alpha-1:\alpha-1:1)$;
  these three coordiantes
  are rationally independent so that $\ell_\infty$ has rationality degree $3$ and therefore it cannot belong
  to any boundary and for Proposition~\ref{p4} it must belong to $\cE$.
\end{proof}
Note that, since $\cE$ is invariant by the $\psi_a$, this way we automatically get a countable
number of explicit elements of $\cE$.
%
\begin{prop}
Let $\{\cD_i\}_{i\in\bN}$ be the set of all stability zones
sorted according to any recursive algorithm, e.g.
$\cD_{(i_1\dots i_n)_3}=\psi_{i_n+1}\circ\dots\circ\psi_{i_1+1}(\cD_{(1:1:1)})$, where
$(i_1\dots i_n)_3$ is the \emph{base 3} expression for the index of the stability zone,
and be $\ell_n$ the label associated to $\cD_n$.
Then $2\log_3^2(1+2n)+1\leqslant \|\ell_n\|^2\leqslant3(1+2n)^{2\log_3\alpha}$ where $\alpha$ is the Tribonacci constant.
\end{prop}
\begin{proof}
  A simple induction argument shows that, at every recursion level $k$, the biggest label belongs,
  modulo permutations of the projective coordinates, to the $k$th stability zone of the steepest descent sequence
  having as first three elements $D_{(0,0,1)}$, $D_{(0,1,0)}$ and $D_{(1,0,0)}$. Since at the $k$th
  recursion level there are $3^k$ stability zones it follows that $\|\ell_{\frac{3^k-1}{2}}\|\leqslant\sqrt{3}\alpha^k$
  and therefore $\|\ell_n\|\leqslant \sqrt{3}(3n)^{\log_3\alpha}$.
  \par
  The slowest sequence which can be formed by picking a stability zone for every recursion level is instead the
  one where $\cD_k$ is the spawn of, say, $D_{(0,0,1)}$, $D_{(0,1,0)}$ and $\cD_{k-1}$. In this case
  indeed $\|\ell_{\frac{3^k-1}{2}}\|=\sqrt{2(k+1)^2+1}$ from which follows immediately the left part
  of the inequality.
\end{proof}
\begin{prop}\label{zoneasympt}
  Be $\cD_\ell$ a stability zone with label $\ell$, $p_\ell$ its
  perimeter and $a_\ell$ its area in $\RPt$.
  Then there exist four constants $A$, $B$, $C$, $D$ such that
  $$
  \frac{A}{\|\ell\|^{\frac{3}{2}}}\leqslant p_\ell
  \leqslant\frac{B}{\|\ell\|}\;,\;\;\frac{C}{\|\ell\|^3}\leqslant a_\ell \leqslant\frac{D}{\|\ell\|^3}
  $$
\end{prop}
\begin{proof}
  The inequalities concerning $p_\ell$ are already known to be true in the general case.
  \par
  To prove the ones relative to the area we change coordinates in $\RPt$ so that the
  three points $(1:0:1)$, $(0:1:1)$, $(1:1:0)$ become $(0:0:1)$, $(1:0:1)$, $(0:1:1)$.
  This way we can use the Lebesgue measure instead of the projective one and we can use
  the fact that $0\leqslant x,y\leqslant z$; now, be $(x_i:y_i:z_i)$, $i=0,1,2$, the projective
  irreducible coordinates of the vertices $A_i$ of the triangle $T$ which encloses $\cD_\ell$,
  so that $l=(\sum_i x_i:\sum_i y_i:\sum_i z_i)$; the areas of $T$ and $\cD_\ell$ are then
  easily computed respectively as $\cA(T)=\Pi_{i=0}^2\frac{1}{z_i}$ and
  $\cA(\cD_\ell)=\Pi_{i=0}^2\frac{1}{z_i+z_{i+1}}$ (where the indices are meant modulo 3).
  \par
  It is convenient for our purposes to consider the maximum norm $\|\ell\|_\infty$ since
  in the region under consideration $0\leqslant x,y\leqslant z$, so that we can assume without
  loss of generality that $x_i,y_i\leqslant z_i$, $i=0,1,2$, and therefore we get trivially that
  $$
  \cA(\cD_\ell)=\frac{1}{(z_0+z_1)(z_1+z_2)(z_2+z_0)}\geq\frac{1}{(z_0+z_1+z_2)^3}=\frac{1}{\|\ell\|^3_\infty}
  $$
  The second part of the inequality comes from the fact that if $\max z_i\leqslant \hbox{mid\ } z_i + \min z_i$
  holds for the first triangle of the algoritm in Prop.~\ref{prop:alg} then it holds for every other
  triangle generated by the algorithm.
  Indeed assume to fix the ideas that $z_0\leqslant z_1\leqslant z_2$; then the new $z$ coordinates of the
  three vertices under the transform $\psi_0(z_0:z_1:z_2)=(z_0:z_1+z_0:z_2+z_0)$ satisfiy trivially
  again the relation
  $z'_2\leqslant z'_0 + z'_1$, and the very same happens for the tranforms under $\psi_1$. In case of
  $\psi_2$ we have $\{z'_0=z_0+z_2,z'_1=z_1+z_2,z'_2=z_2\}$ so that the new larger $z$ is now $z'_1$
  but again $z'_1\leqslant z'_2+z'_0$.
  Then finally
  $$\frac{\|\ell\|^3_\infty}{(z_0+z_1)(z_1+z_2)(z_2+z_0)}=
  \Bigl(1+\frac{z_0}{z_1+z_2}\Bigr)\Bigl(1+\frac{z_1}{z_2+z_0}\Bigr)
  \Bigl(1+\frac{z_2}{z_0+z_1}\Bigr)\leqslant8$$
  so that $\cA(\cD_\ell)\leqslant8/\|\ell\|^3_\infty$. The final statement now follows from the fact that all norms
  are equivalent in finite dimension.
\end{proof}
\begin{figure}
  \begin{center}
    \begin{tabular}{cc}
      \includegraphics[width=5.85cm]{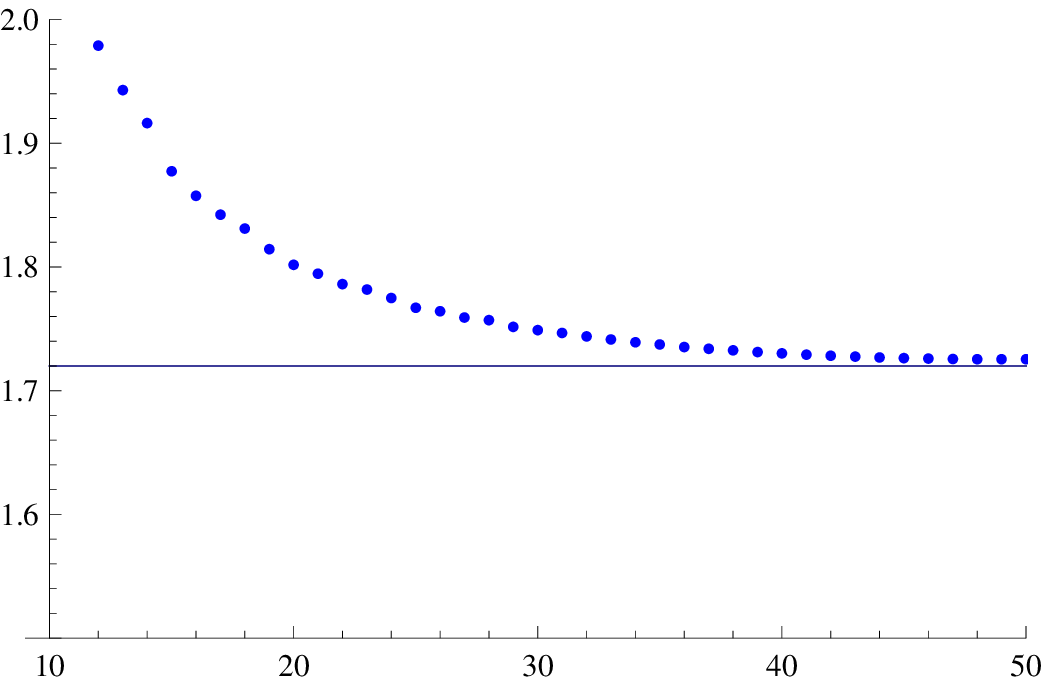}&
      \includegraphics[width=5.85cm]{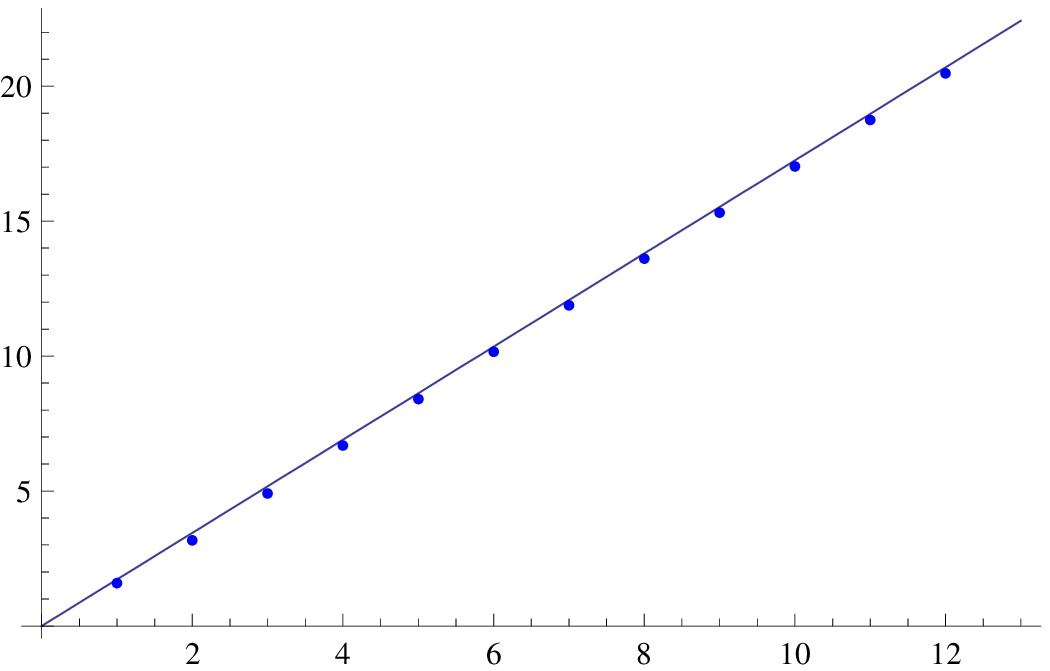}\\
    \end{tabular}
  \end{center}
  \caption{%
    \small
    (a) Plot of $d_n=2-\log_{1.2}V(\cE_{1.2^{-n}})$ against $n$, where $V(\cE_{1.2^{-n}})$ is the volume
    of the neighborhood of $\cE$ of radius $1.2^{-n}$; for ``well-behaving'' fractals this sequence
    converges to their Minkowsky dimension. The horizontal line shown in the picture has height $1.72$.
    (b) Plot of $\log_2 N_n$ against $n$, where $N_n$ is the number of squares of side $2^{-n}$ which
    are not completely inside one of the 797161 stability zones found at the 12th recursion level. The angular
    coefficient of the straight line fitting the data in the picture is again $1.72$.
  }
  \label{fig:numfracdim}
\end{figure}
We could not find any way to evaluate exact non-trivial bounds for the fractal dimension
of $\cE$ but numerical calculations suggest that the dimension be smaller than 1.8.
\par
A simple way to evaluate numerically fratal dimensions is using the \emph{Minkowsky dimension},
namely the limit
$$
\dim_M\cE=2-\lim_{\epsilon\to0}\frac{\log V(\cE_\epsilon)}{\log \epsilon}
$$
where $V(\cE_\epsilon)$ is the surface of the $\epsilon$ neighborhood of $\cE$.
In order to do that we use the fact that, if $A$ is the area of $\Delta$ and $p$ its
perimeter, then
$$
V(\cE_\epsilon) = p\epsilon + \epsilon\sum_{i=1}^{k_\epsilon} p_i + A - \sum_{i=1}^{k_\epsilon} a_i +
\epsilon^2(\pi-\sum_{i=1}^{k_\epsilon}\frac{p_i^2}{4a_i})
$$
where $k_\epsilon$ is the integer such that $\rho_{k_\epsilon+1}\leqslant\epsilon\leqslant\rho_{k_\epsilon}$
and $\rho_k$ is the radius of the inscribed circle to the triangle $\cD_k$.
In order to avoid infinities we make a projective change of coordinates so that the triangle
$\Delta$ has vertices in $(0:0:1)$, $(1:0:1)$ and $(0:1:1)$. In Figure~\ref{fig:numfracdim}(a)
we show the numerical results we got by evaluating the volume of the neighborhoods of
$\cE$ of radii $r_n=1.2^{-n}$ for $n=1,\ldots,50$, which suggests a fractal dimension
between $1.7$ and $1.8$.
\par
A second simple way is to evaluate an upper bound for the \emph{box-counting dimension},
namely the limit
$$
\dim_B\cE=\lim_{\epsilon\to0}\frac{\log N_\epsilon(\cE)}{-\log \epsilon}
$$
where $N_\epsilon(\cE)$ is the smallest number of squares of side $\epsilon$ needed to
cover $\cE$ (even in this case we make the same change of coordinates to avoid infinities).
In Figure~\ref{fig:numfracdim}(b) we show the results relative to covering with squares
of side $2^{-n}$, $n=1,\cdots,12$, the complement of the 797161 triangles obtained by
applying 12 times the recursion algorithm. Again we get a clear indication of the fractal
dimension to be between 1.7 and 1.8.
\section{Numerical generation of stability zones}
\begin{figure}
  \begin{center}
    \begin{tabular}{cc}
      \includegraphics[width=5.85cm]{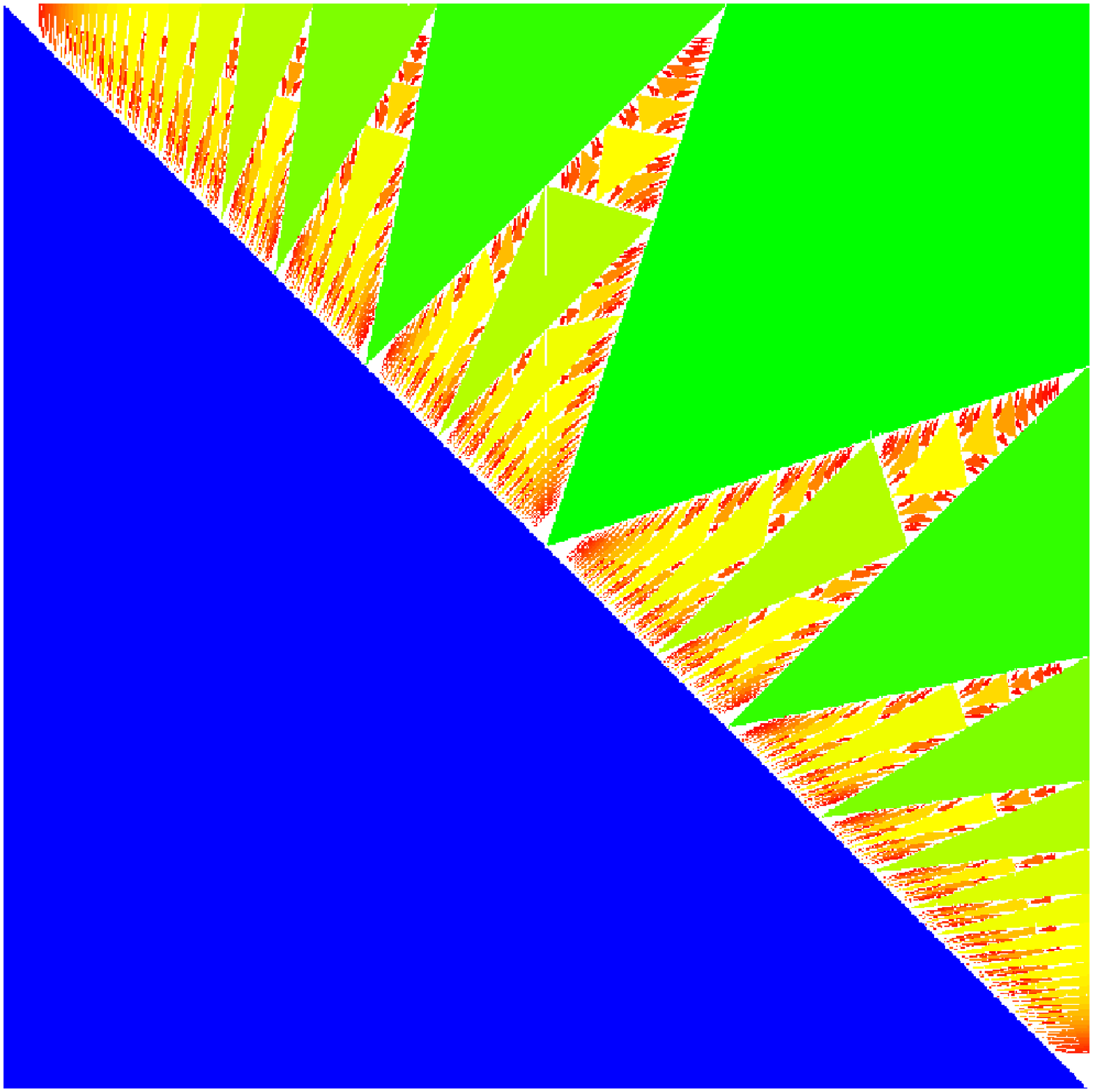}&
      \includegraphics[width=5.5cm,viewport=0 128 566 667,clip]{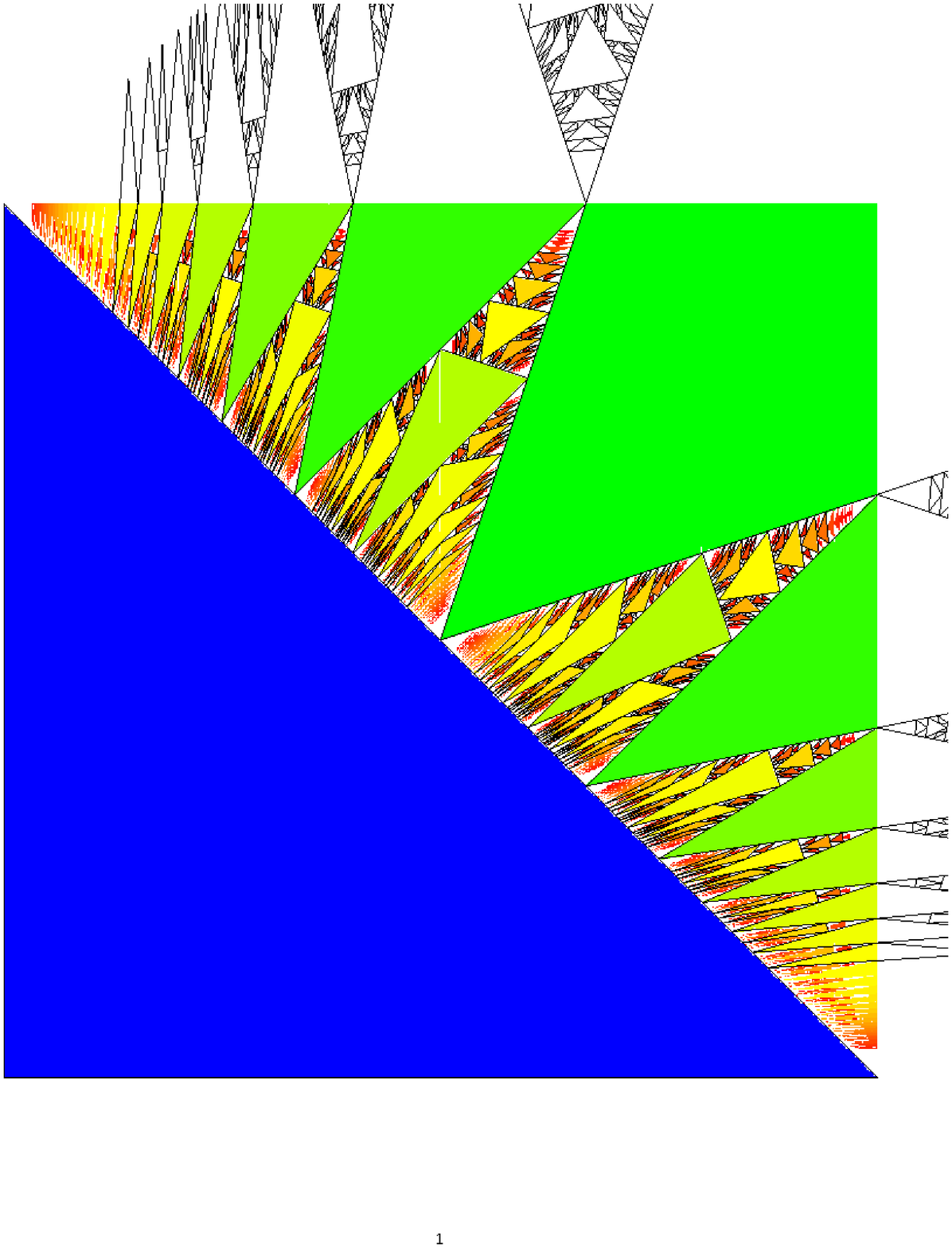}\\
    \end{tabular}
  \end{center}
  \caption{%
    \small
    (a) Picture of the fractal in the square $[0,1]^2$ in the $z=1$ chart of $\RPt$ obtained
    by evaluating the soul of every stability zone in the lattice $(n,m,N)$, $n,m=1,\ldots,10^3$, $N=10^3$;
    (b) the same picture compared with the analytical boundaries of the $797161$ stability zones of the first
    $12$ levels of recursion. Out of the 961367 rational directions for which a label was numerically
    found by the algorithm, 455654 belong to those $797161$ stability zones and all of them turn out to lie
    within the analytical boundaries of the corresponding zone.
  }
  \label{fig:numfrac}
\end{figure}
As an important byproduct of the study of stability zones
of $\cC$, we were able for the first time
to compare very accurately the results of the NTC code~\cite{NTC},
used in~\cite{DeL03a,DeL03b,DeL06}
to generate approximations of surfaces' generality stability zones against their
analytical boundaries.
\par
Indeed no simple algorithm to generate the analytical equations of the boundaries of stability zones
for a generic function is known so far but we know that all directions
belonging to the same stability zone share the same soul and therefore it is possible
to get an approximate picture of the set of generalized stability zones by evaluating the soul
of some (big) set of rational directions. For example in all cases examined so far, thanks
to the high level of symmetry, we could limit our analysis to the directions contained in the
the triangle with sides $(0:0:1)$, $(1:1:1)$ and $(1:0:1)$ (concretely, to all directions
$(m:n:N)$, $0<n<m<N$, for $N=10^2,10^3,10^4$).
\par
Note that rational directions in this setting are of paramount importance because
their (non-critical) leaves are compact (in $\Tt$) and therefore can be in principle
approximated with error as small as wished and therefore the corresponding soul can
be, at will, evaluated \emph{exactly} through numerical calculations.
Moreover, rational directions are dense in every stability zone and therefore (in principle)
we do not lose any picture detail by restricting our analysis to them.
\par
Below we present the algorithm we use to retrieve the soul (if any) of a rational direction
$(m:n:N)$ in this particular case, where we have all saddles of ``monkey'' type:
\begin{itemize}
\item[{\textbf N0}] choose a representative $b_i$, $i=1,2,3$, for the cycles
  on $\cC$ which are respectively homologous in $\Tt$ to the $x$, $y$ and $z$ axes;
\item[{\textbf N1}] retrieve the intersection between $\cC$ and the plane $$m(x-1/4)+n(y-1/4)+N(z-1/4)=0\,;$$
\item[{\textbf N2}] follow the three critical loops and, if no saddle connection is detected,
  store them in the variables $c_{1,2,3}$, otherwise exit;
\item[{\textbf N3}] evaluate the homology class of $c_{1,2,3}$ in $\Tt$;
\item[{\textbf N4}] if exactly one of the three loops is non-zero in $\Tt$ then evaluate its
  intersection numbers with the $\{b_i\}$ and set this triple as the soul corresponding to
  the direction $(m,n,N)$, otherwise exit.
\end{itemize}
\par
The result of sampling the triangle $\cT$ with a $10^3 \times 10^3$ lattice are shown in
Figure~\ref{fig:numfrac} and turn out to be in perfect agreement with the analytical boundaries.
An evaluation of the fractal box dimension based on these numerical data leads to a value of
about $1.7$,
which is also in very good agreement with the evaluation
obtained from of the analytical boundaries made in sec.~\ref{sec:asympt}.
%
%
%
\section{Acknowledgments}
The authors gladly thanks the IPST (www.ipst.umd.edu) and the Dept. of Mathematics of the UMD (USA)
(www.math.umd.edu) for their hospitality in the Spring Semester 2007 and for financial support.
Numerical calculation were made on Linux PCs kindly provided by the UMD Mathematics Dept.
and by the Cagliari section of INFN (www.ca.infn.it) which also provided financial support
tho the first author. The authors also warmly thank S.P. Novikov and B. Hunt for several
fruitful discussion during their stay at UMD.
\bibliography{dc}
\end{document}